\documentclass{amsart}
\usepackage{amsmath,amssymb,amscd,latexsym,enumitem,amsthm, hyperref, float}

\usepackage[all]{xy}
\usepackage{color}

\newcommand{\N}{{\mathbb{N}}}
\newcommand{\R}{{\mathbb{R}}}

\DeclareMathOperator{\Q}{\mathcal{Q}}

\DeclareMathOperator{\id}{\textup{id}}
\DeclareMathOperator{\into}{\hookrightarrow}

\DeclareMathOperator{\Ell}{Ell}

\newcounter{number}[section]

\newenvironment{nummer}{\refstepcounter{number}{\noindent\arabic{section}.\arabic{number}}}{}

\newcommand{\bn}{\noindent \begin{nummer} \rm}
\newcommand{\en}{\end{nummer}}

\newenvironment{ntheorem}{\noindent {\sc  Theorem:} \it}{}
\newenvironment{nlemma}{\noindent {\sc  Lemma:} \it}{}
\newenvironment{nprop}{\noindent {\sc  Proposition:} \it}{}
\newenvironment{ndefn}{\noindent {\sc  Definition:} \it}{}
\newenvironment{ncor}{\noindent {\sc  Corollary:} \it}{}
\newenvironment{nconj}{\noindent {\sc Conjecture:} \it}{}
\newenvironment{nremark}{\noindent {\sc  Remark:}}{}
\newenvironment{nremarks}{\noindent {\sc  Remarks:}}{}

\begin{document}

\title[Minimal dynamics and the Jiang--Su algebra]{Constructing minimal homeomorphisms on point-like spaces and a dynamical presentation of the Jiang--Su algebra}
\author[R.J. Deeley, I.F. Putnam, K.R. Strung]
{Robin J. Deeley \and
Ian F. Putnam \and
Karen R. Strung}
\address{Universit\'e Blaise Pascal, Clermont-Ferrand II, Laboratoire de Math\'ematiques, Campus des C\'ezeaux 63177 Aubi\`ere cedex, France} 
\email{robin.deeley@gmail.com}
\address{Department of Mathematics and Statistics,
University of Victoria,
Victoria, B.C., Canada V8W 3R4} 
\email{ifputnam@uvic.ca}
\address{Instytut Matematyczny Polskiej Akademii Nauk, ul. \'{S}niadeckich 8, 00-656 Warszawa, Poland} 
\email{kstrung@impan.pl}

\date{\today}
\thanks{The first author was supported by ANR Project SingStar. The second author was supported by an NSERC
 Discovery Grant.}

\setcounter{section}{-1}

\begin{abstract}
The principal aim of this paper is to  give a dynamical presentation of the Jiang--Su algebra. Originally constructed as an inductive limit of prime dimension drop algebras, the Jiang--Su algebra has gone from being a poorly understood oddity to having a prominent positive  role in George Elliott's classification programme for separable, nuclear $\mathrm{C}^{*}$-algebras. Here, we exhibit an \'{e}tale equivalence relation whose groupoid $\mathrm{C}^{*}$-algebra is isomorphic to the Jiang--Su algebra. The main ingredient is the construction of  minimal homeomorphisms on  infinite, compact metric spaces, each having  the same cohomology as a point. This
 construction is also of interest in dynamical systems. 
Any self-map of an infinite, compact space with the same cohomology 
as a point has  Lefschetz number one. Thus, 
if such 
a space were also
to satisfy some regularity hypothesis 
(which our examples do not), 
then  the Lefschetz--Hopf Theorem would imply that it 
does not admit a  minimal homeomorphism.
\end{abstract}

\maketitle

\section{Introduction}
The fields of operator algebras and dynamical systems have a long history of mutual influence. On the one hand, dynamical systems provide interesting examples of operator algebras and have often provided techniques which are successfully imported into the operator algebra framework. On the other hand, results in operator algebras are often of interest to those in dynamical systems. In ideal situations, significant information is retained when passing from dynamics to operator algebras, and vice versa.

This relationship has been particularly interesting for the classification of $\mathrm{C}^*$-algebras. An extraordinary result in this setting is the  classification, up to strong orbit equivalence, of the minimal dynamical systems on a Cantor set and the corresponding $K$-theoretical classification of the associated crossed product $\mathrm{C}^*$-algebras \cite{GioPutSkau:orbit, Put:TopStabCantor}.   Classification for separable, simple, nuclear $\mathrm{C}^*$-algebras remains an interesting open problem. To every simple separable nuclear $\mathrm{C}^*$-algebra one assigns  a computable set of invariants involving $K$-theory, tracial state spaces, and the pairing  between these objects.  George Elliott conjectured that for all such $\mathrm{C}^*$-algebras, an isomorphism at the level of invariants, now known as Elliott invariants, might be lifted to a $^*$-isomorphism at the level of $\mathrm{C}^*$-algebras.
 
A remarkable number of positive results have been obtained. However, examples---including examples of some crossed product $\mathrm{C}^*$-algebras arising from minimal dynamical systems---have also shown pathologies undetectable by the Elliott invariant \cite{GioKerr:Subshifts, Ror:simple, Toms:classproblem, Vil:perforation}. One such algebra, constructed by Xinhui Jiang and Hongbing Su,  gives a $\mathrm{C}^*$-algebra $\mathcal{Z}$ with invariant isomorphic to that of $\mathbb{C}$ \cite{JiaSu:Z}. The importance of the Jiang--Su algebra for the classification programme cannot be understated:  in the case that a $\mathrm{C}^*$-algebra $A$ has weakly unperforated $K_0$-group (for a definition, see for example \cite[Definition 6.7.1]{Bla:newk-theory}), its Elliott invariant is isomorphic to the Elliott invariant of $A \otimes \mathcal{Z}$.  The original Elliott conjecture then predicts that for any simple separable nuclear $\mathrm{C}^*$-algebra with weakly unperforated $K_0$-group, $A \cong A \otimes \mathcal{Z}$. In such a case $A$ is said to be a $\mathcal{Z}$-stable $\mathrm{C}^*$-algebra. So far, each counterexample to Elliott's conjecture involves two $\mathrm{C}^*$-algebras, one of which is not $\mathcal{Z}$-stable. This leads to the following revised conjecture:

\bn
\begin{nconj}[Toms--Winter, 2008]
Let $A$ and $B$ be simple separable unital nuclear $\mathrm{C}^*$-algebras. Suppose that $A$ and $B$ are $\mathcal{Z}$-stable and have isomorphic Elliott invariants. Then $A \cong B$. 
\end{nconj}
\en

More recently, there has been significant interest in transferring $\mathrm{C}^*$-algebraic regularity properties to the language of topological dynamics with the aim of showing that the appropriate properties pass from dynamical system to  associated crossed product $\mathrm{C}^*$-algebra. Much of the motivation for this has come from the theory of von Neumann algebras, which has close ties to ergodic theory. Currently, the classification programme for $\mathrm{C}^*$-algebras is seeing rapid advancement by adapting results for von Neumann algebras to the setting of $\mathrm{C}^*$-algebras. For a good discussion on this interplay, we refer the reader to \cite{SatWhiWin:NucDimZ}. Of interest here are the comparisons between $\mathcal{Z}$ and its von Neumann counterpart, the hyperfinite $\textup{II}_1$-factor, $\mathcal{R}$.

Like the Jiang--Su algebra, $\mathcal{R}$ is strongly self-absorbing (see \cite{TomsWin:ssa}), absorbed (after taking tensor products) by certain factors with a particularly nice structure, and can be characterised uniquely in various abstract ways. Francis Murray and John von Neumann realised $\mathcal{R}$ using their group measure space construction, the von Neumann algebra version of the crossed product of a commutative $\mathrm{C}^*$-algebra by the integers \cite{MvN:Rings4}. Its ties to ergodic theory were deepened in \cite{Con:class}, where Alain Connes proves that $\mathcal{R}$ can be realised dynamically by any measure space $(X, \mu)$ with a probability measure preserving action of a discrete amenable group $G$. Such a dynamical presentation for $\mathcal{Z}$ has so far been missing from the $\mathrm{C}^*$-algebraic theory. 

In light of this, it has become increasingly important to find a suitable dynamical interpretation of $\mathcal{Z}$. In this paper, we construct such a presentation of the Jiang--Su algebra via a minimal \'{e}tale equivalence relation (that is, an equivalence relation with countable dense equivalence classes.) See \cite{Ren:groupoid} for more about these groupoids. Along the way, we meet an old question in dynamical systems:  which compact, metric spaces admit minimal homeomorphisms?

Of course, many well-known systems provide positive answers (Cantor sets from odometers, the circle from irrational rotations). Perhaps the most famous positive result is that of Albert Fathi and Michael Herman who exhibited minimal, uniquely ergodic
 homeomorphisms on all odd-dimensional spheres of dimension greater than one \cite{FatHer:Diffeo} and Alistair Windsor's subsequent generalization to arbitrary numbers of ergodic measures \cite{Wind:not_uniquely_ergo}. In fact, we will use these results in a crucial way in our construction.

There are also negative results. Perhaps the most famous,  and the most relevant for our discussion, is the Lefschetz--Hopf theorem (see for example \cite{Brown:LefFixPt}) which asserts that for ``nice'' spaces (for example, absolute neighbourhood retracts), the cohomology of a space  may contain enough information  to conclude  that any continuous self-map of the  space has a periodic  point. If we also ask that the space be infinite, then it does not admit a minimal homeomorphism. An example where this holds is any even-dimensional sphere. The same conclusion holds for  any contractible absolute neighbourhood retract (ANR).
More generally, it also applies to any ANR  whose cohomology is the same as a point. 

Here we build minimal homeomorphisms $\zeta$ on ``point-like'' spaces:  infinite, compact metric spaces with the same cohomology and $K$-theory as a point. In fact, our spaces are inverse limit of ANR's so while our results are positive, they sit perilously close to the Lefschetz--Hopf trap.  For a survey on fixed point properties, see \cite{Bing:Fixed}.

We construct such minimal dynamical systems with any prescribed number of ergodic Borel probability measures. When the system $(Z, \zeta)$ is uniquely ergodic, the resulting $\mathrm{C}^*$-algebra  $C(Z) \rtimes_{\zeta} \mathbb{Z}$ then has the same invariant as $\mathcal{Z}$, except for its nontrivial $K_1$-group. For such a crossed product, nontrivial $K_1$ is unavoidable (the class of the unitary implementing the action is nontrivial). However, upon breaking the orbit equivalence relation across a single point, $K_1$ disappears while the rest of the invariant remains the same. Now the $\mathrm{C}^*$-algebra arising from this equivalence relation does in fact have the correct invariant and we are able to use classification theory to conclude that it must be $\mathcal{Z}$. 

Our construction itself is more general, and in fact we are able to produce $\mathrm{C}^*$-algebras isomorphic to any simple inductive limit of prime dimension drop algebras with an arbitrary number of extreme tracial states, as constructed in \cite{JiaSu:Z}. From the $\mathrm{C}^*$-algebraic perspective such a construction is interesting: even the range of the invariant for such $\mathrm{C}^*$-algebras remains unknown. In the uniquely ergodic case, classification for the resulting crossed product follows from \cite{TomsWinter:PNAS, TomsWinter:minhom} (which uses the main result in \cite{StrWin:Z-stab_min_dyn}). Without assuming unique ergodicity, we may appeal to Lin's generalisation \cite{Lin:Spheres} of the third author's classification result for products with Cantor systems \cite{Str:XxSn}, to show all our minimal dynamical systems result in classifiable crossed products. We note that as we were finishing this paper,  Lin posted a classification theorem for all crossed product $\mathrm{C}^*$-algebras associated to minimal dynamical systems with mean dimension zero \cite{Lin:MinDyn}, but our results do not rely on his proof.

In Section~\ref{construct}, we start with  a minimal diffeomorphism, $\varphi$, on an $d$-sphere, for odd $d>1$, 
which we denote by $S^{d}$.  This is a logical place 
to begin since the cohomology of the sphere differs from 
that of a point only in dimension $d$. 
From this, we construct a space, $Z$, together
 with a minimal homeomorphism, 
$\zeta$. This system is an extension of $(S^{d}, \varphi)$; that is, there
a a factor map from $(Z, \zeta)$ onto $(S^{d}, \varphi)$.
 The space $Z$ is an infinite compact finite-dimensional metric space with the same cohomology and $K$-theory as a single point. From our minimal dynamical system we show in Section~\ref{Jiang-Su} that the $\mathrm{C}^*$-algebra of the associated orbit-breaking equivalence relation is isomorphic to the Jiang--Su algebra, assuming we have begun with uniquely ergodic $(S^d, \varphi)$.  In Section~\ref{classification}, we show that, with any number of ergodic measures, the associated transformation group $\mathrm{C}^*$-algebras and their orbit-breaking subalgebras can be distinguished by their tracial states spaces and are all isomorphic to direct limits of dimension-drop algebras. Finally, in Section~\ref{outlook}, we make some comments on further questions.

\section{Constructing the system} \label{construct}

In this section, we fix $d > 1$ odd and a minimal diffeomorphism $\varphi: S^d \to S^d$, but we remind the reader that the space constructed does in fact depend on which minimal dynamical system $(S^d, \varphi)$ we use. In particular, we may choose $(S^d, \varphi)$ to have any number of ergodic probability measures \cite{Wind:not_uniquely_ergo}. We fix an orientation on $S^d$ and note that $\varphi$ is orientation-preserving; otherwise the system would have a fixed point.  

\bn
  \begin{nlemma}
\label{Z:3}
Let $x$ be any point of $S^{d}$ and 
$v$ any non-zero tangent vector at $x$.
There exists 
\[
\lambda : [0, 1 ] 
 \rightarrow S^{d}
\] 
satisfying the following:
\begin{enumerate}
\item
$\lambda(0) = x, \lambda'(0) = v, \lambda(1) = \varphi(x), 
\lambda'(1) = D\varphi(v)$,
\item $\lambda$ is a $C^{1}$-embedding, in particular $\lambda'(t) \neq 0$, for all $ t $ in $[0,1]$,
\item For all $n \neq 0$,
\[
\varphi^{n}( \lambda ([0, 1 )  ) 
\cap \lambda ([0, 1 ) ) = \emptyset.
\]
\end{enumerate}
\end{nlemma}

\begin{proof}
The space of all $C^{1}$-maps from $[0,1]$ into $S^{d}$ which satisfy the first condition is a non-empty, complete metric space with the metric from the $C^{1}$-norm.  Let $\Lambda$ be the subset of these also
  satisfying the  second condition of the conclusion. This is clearly non-empty and open (see for example \cite[Chapter 2, Theorem 1.4]{Hirsch:DiffTop}). Thus, $\Lambda$ is a Baire space. 

We need to establish the existence of a map satisfying the third
condition as well. 
We treat the end points of $[0,1]$ 
separately.

For each  integer $n \neq 0, 1$, 
let $\Lambda^{0}_{n}$ be those elements of 
$\Lambda$ such that  
$ \varphi^{n}(x) \notin \lambda[0, 1 ]$. 
This is clearly an open dense 
subset of $\Lambda$. In particular, 
 the intersection
over all $n \neq 0, 1$, which we denote 
$\Lambda^{0}_{\infty}$, is a dense $G_{\delta}$ in
$\Lambda$, and hence also a Baire space.

Fix $\lambda$ in $\Lambda$.
 For each $n \geq 1, k \geq 2$, define 
\begin{eqnarray*}
R_{n,k}(\lambda) & =  & \{ ( s_{1}, s_{2}, \ldots, s_{k}) \mid s_{i}
\in [0, 1], \varphi^{n}(\lambda(s_{i}) ) = \lambda(s_{i+1}), 
1 \leq  i < k \}, \\
X_{n,k}(\lambda) & =  & \{ s_{1} \mid  ( s_{1}, s_{2}, \ldots, s_{k}) \in R_{n,k}(\lambda) \}.
\end{eqnarray*}

Let us start with some simple observations. 
\begin{enumerate}
\item  $X_{n,k}(\lambda) \supseteq X_{n,k+1}(\lambda)$, for all  $n \geq 1, k \geq 2$. 
\item It follows from the first condition that
$(0,1) \in R_{1,2}(\lambda)$ and that it is an isolated point; 
$0 \in X_{1,2}(\lambda)$ and is an isolated point.
\item It follows from the preceding two
 properties that, for any $k \geq 2$,
if $0 $ is in $ X_{1,k}(\lambda)$, then it is an isolated point.
 \item For all $n, k$,  
 $R_{n,k}(\lambda)$ is closed in $[0, 1]^{k}$ and $X_{n,k}(\lambda)$ is closed in $[0,1]$. 
 \item For all $n \geq 1, k \geq 2$, we
  have 
  \[
  X_{n,k}(\lambda) = \emptyset \Leftrightarrow 
R_{n,k}(\lambda) = \emptyset \Rightarrow 
R_{n,k+1}(\lambda) = \emptyset \Leftrightarrow 
 X_{n,k+1}(\lambda) = \emptyset.
 \]
 \item   
 $\{ \lambda \in \Lambda \mid X_{1,2}(\lambda) = \{ 0 \} \}$ is open in 
 $\Lambda$.
 \item  
  For any  $n \geq 1, k \geq 2$, 
 $\{ \lambda \in \Lambda \mid X_{n,k}(\lambda) = \emptyset \}$
  is open in 
 $\Lambda$. 
 \item If $X_{1,2}(\lambda) = \{ 0\}$ and 
 $X_{n,2}(\lambda) = \emptyset$,
 for all $n \neq 0, 1$, then $\lambda$ 
 satisfies the last condition of the conclusion of the lemma.
 \end{enumerate}

 Our first important claim is that, for any $n \geq 1$, 
 there exists  $k > 1$ with  $R_{n,k}(\lambda) = \emptyset$.
  Suppose the contrary. We note there is 
 an obvious  map 
 from $R_{n,k}$ to $R_{n,k-1}$ and we may form the 
 inductive limit of this  system. Using the compactness of 
 $R_{n,k}$, we  see that if each $R_{n,k}$ is non-empty, we may 
 find a sequence $s_{1}, s_{2}, \ldots$ such that 
 $\varphi^{n}(\lambda(s_{i})) = \lambda(s_{i+1})$, for 
 all $ i \geq 1$.  But this means that the forward orbit 
 of $\lambda(s_{1})$ under $\varphi$ is contained 
 in 
 $\cup_{i=0}^{n-1} \varphi^{i}(\lambda[0,1])$, which 
 is a closed subset of $S^{d}$. It is non-empty and 
 cannot be all of $S^{d}$ 
 on dimensionality grounds. This then contradicts the 
 minimality of $\varphi$.
 
Our second claim is the following. Suppose
that  $\lambda$ is in $\Lambda$, 
$n \geq 1, k \geq 2, (n,k) \neq (1,2)$
 are such that $X_{n, k+1}(\lambda) = \emptyset$.
Then for any $\epsilon > 0$, there is 
$\mu$ in $\Lambda$ with $\Vert \lambda - \mu \Vert_{1} < \epsilon $ 
and $X_{n, k}(\mu) = \emptyset$.

In view of the first claim above and 
the fact that $X_{n, k+1}(\lambda) = \emptyset$ is 
an open property, there 
is no  loss of generality if we assume that $\lambda$
is in $\Lambda^{0}_{\infty}$. The immediate consequence of this 
is that $0,1$ are not in $X_{n,k}(\lambda)$.

Suppose that  $s$ is in $[0,1]$ with $\varphi^{n}(\lambda(s))$  in 
$\lambda(X_{n,k}(\lambda)))$, then $\varphi^{n}(\lambda(s)) = s_{1}$
with $(s_{1}, s_{2}, \ldots, s_{k})$ in $R_{n,k}(\lambda)$ and it follows
that $(s, s_{1},s_{2}, \ldots, s_{k})$ is in $R_{n,k+1}(\lambda)$ 
which contradicts our hypothesis that  $X_{n, k+1}(\lambda) = \emptyset$.
We conclude that 
the sets $\lambda(X_{n,k}(\lambda)))$ and $\varphi^{n}(\lambda[0,1])$
 must be disjoint. Without loss of generality assume that $\epsilon$ is 
 strictly less than half the distance between these two compact sets.
 
 For each $s$ in $X_{n,k}(\lambda)$, select $0 < a_{s} < s < b_{s} < 1$
 such that $\lambda(a_{s}, b_{s})$ is contained in the ball 
 of radius $\epsilon/2$ about $\lambda(s)$. These open intervals 
 cover $X_{n, k}(\lambda)$. We may extract a finite subcover. If these 
 intervals overlap, we may replace them with their unions to 
 obtain $0 < a_{1} < b_{1} < a_{2} < \cdots < b_{n} < 1$ with 
 the union of the $(a_{i}, b_{i})$, which we denote by $U$,
  covering $X_{n,k}(\lambda)$. Observe that 
 this means the points $a_{i}, b_{i}$ are not in   $X_{n,k}(\lambda)$.

 Based on dimensionality, we may make an 
 arbitrarily small  $C^{1}$-perturbation 
 of  $\lambda$ on  $U$, which we call $\mu$, 
 not changing the value
 or derivative at the endpoints, so that 
 the image is disjoint from 
 $\varphi^{-n} \circ \lambda([0, 1] - U)$.
  To be slightly more precise, the $\mu$ can be chosen 
  from $\Lambda$ so that 
 $\Vert \lambda - \mu \Vert_{1} < \epsilon/2 $ and 
 $\Vert \varphi^{n}  \circ \lambda -  \varphi^{n}  \circ \mu \Vert_{1} < \epsilon $.

  We claim this 
 $\mu$ satisfies  $X_{n, k}(\mu) = \emptyset$.
 Suppose to the contrary that 
 $(s_{1}, s_{2},  \ldots, s_{k})$ is in $X_{n, k}(\mu)$.
 If, for some $j < k$, $s_{j}$ is in $U$, then
 $\mu(s_{j})$ is not in $\varphi^{-n} \circ \lambda([0,1]-U)$
 or equivalently, $\varphi^{n}(\mu(s_{j}))$ is not in 
 $\lambda([0,1] -U) =  \mu([0,1] -U)$. On the other hand, $\varphi^{n}(\mu(s_{j}))$
 is within $\epsilon$ of $\varphi^{n}(\lambda[0,1])$
 and hence outside the ball of radius $\epsilon$ of 
 $\lambda(X_{n,k}(\lambda))$, which contains
 $\mu(U)$. Between the two cases, we have shown that 
  if $s_{j}$ is in $U$, then $\varphi^{n}(\mu(s_{j}))$ 
  is not in $\mu([0,1])$. This contradicts 
  $\varphi^{n}(\mu(s_{j}) ) = \mu(s_{j+1})$.
  
 The only remaining case is that $s_{1}, s_{2}, \ldots, s_{k-1}$
 all lie in $[0,1] - U$. But on this set, $\mu = \lambda$
 and so we have $(s_{1}, s_{2},  \ldots, s_{k})$ 
 is in $X_{n, k}(\mu) = X_{n,k}(\lambda)$ and hence $s_{1}$ 
 is in $X_{n,k}(\lambda) \subseteq U$, 
 a contradiction.
 
 Our third claim is a minor variation of the 
 second to deal with the case $n=1, k=2$, 
 since $0$ always lies in $X_{1,2}(\lambda)$.
 Suppose
that  $\lambda$ is in $\Lambda$ 
 such
 that $X_{1, 3}(\lambda) = \emptyset$.
 Then for any $\epsilon > 0$, there is 
$\mu$ in $\Lambda$ with $\Vert \lambda - \mu \Vert_{1} < \epsilon $ 
and $X_{1,2}(\mu) = \{ 0 \}$. The idea is that 
$0$ will be an isolated point of 
$X_{1,2}(\mu) $ and we can simply repeat the rest 
of the argument above replacing $X_{n, k}(\mu) $
by $X_{1,2}(\mu) - \{ 0 \}$.

 Our fourth claim is that, for any $n \neq 0, 1$,
 $\{ \lambda \mid X_{n,2}(\lambda) = \emptyset \}$ is dense 
 in $\Lambda$. Let $\lambda$ be in $\Lambda$ and let 
 $\epsilon > 0$. From our first claim, we may find 
 $k$ with  $X_{n,k}(\lambda) = \emptyset $. Next, use the
 second claim to find $\lambda_{1}$ within $\epsilon /2$ 
 of $\lambda$ with $X_{n,k-1}(\lambda_{1}) = \emptyset $. 
 Apply the second claim again to find $\lambda_{2}$ within 
 $\epsilon/4$ of $\lambda_{1}$  
 with $X_{n,k-2}(\lambda_{2}) = \emptyset $. 
 Continuing in this way, we will end up with 
 $X_{n,2}(\lambda_{k-2}) = \emptyset $ and 
 \[
 \Vert \lambda - \lambda_{k-2} \Vert_{1} \leq 
 \frac{\epsilon}{2} +   
  \frac{\epsilon}{4} +  \ldots +
 \frac{\epsilon}{2^{k-2}}    < \epsilon.
\]

The fifth claim is a minor variation of the third: 
$\{ \lambda \mid X_{1,2}(\lambda) = \{ 0 \} \}$ is dense 
 in $\Lambda$. The proof is the same as that of the fourth claim, 
 using the third claim in place of the second.
 
 As a result of all of this, the set 
 \[
 \{ \lambda \mid X_{1,2}(\lambda) =
  \{ 0 \} \} \cap \left( \cap_{n \neq 0,1} \{ \lambda \mid 
 X_{n,2}(\lambda) = \emptyset \} \right)
 \]
 is dense in $\Lambda$ and hence non-empty. These maps satisfy all the desired conditions.  \end{proof}  \en
 
Given a map $\lambda : [0,1] \to S^d$ satisfying Lemma~\ref{Z:3}, we define a second map 
\[ \lambda_{\mathbb{R}} : \mathbb{R} \to S^d,  \quad \lambda_{\R}(s) = \varphi^{\lfloor{s}\rfloor}(\lambda(s \textup{ mod } 1)) \]
where $\lfloor s \rfloor$ denotes the floor function applied to the real number $s$. The conditions in Lemma~\ref{Z:3} ensure $\lambda_{\R}$ is also a $C^1$-embedding.

\bn  \begin{nlemma}
Let $F: [1,2] \times \R^{d-1} \rightarrow \R^d$ be an orientation preserving $C^1$-embedding such that, for each $s \in [1,2]$,
$$F(s,0)=(s,0).$$
Then, there exists $G:[-1,2]\times \R^{d-1} \rightarrow \R^d$ a continuous function satisfying the following
\begin{enumerate}
\item for $s \in [-1, 2]$, $G(s,0)=(s,0)$;
\item for $x \in \R^{d-1}$ and $s \in [-1, 0]$, $G(s, x)= (s , x)$;
\item for $x \in \R^{d-1}$ and $s \in [1, 2]$, $G(s, x) = F(s,x)$;
\item there exists $\delta>0$ such that $G|_{[-1, 2] \times \overline{B^{d-1}(\delta)}}$ is injective.
\end{enumerate}
\end{nlemma}

\begin{proof}
Note that we have
$$F(s, x)= (s+ F_1(s, x), F_2(s,x))$$
where $F_1(s,0)=0$ and $F_2(s,0)=0$. Observe that 
\[ DF(s,x) =  \left( \begin{array}{cc} 1 + \frac{\partial F_1}{\partial s}(s,x) & \frac{\partial F_1}{\partial x}(s,x) \\
\frac{\partial F_2}{\partial s}(s,x) & [DF_2(s,x)]_{i, j = 2}^{d-1}(s,x) \end{array} \right) \]
where by abuse of notation $\frac{\partial F_1}{\partial x}$ denotes $(\frac{\partial F_1}{\partial x_1}, \dots, \frac{\partial F_1}{\partial x_{d-1}})$ where $x= (x_1, \dots, x_{d-1})\in \R^{d-1}$ and $[DF_2(s,x)]_{i j = 2}^{d-1}$ denotes the minor of the Jacobian matrix of $F_2(s,x)$ obtained by removing the first row and column. The conditions $F_1(s,0)=0$, respectively $F_2(s,0)=0$ imply that $\frac{\partial F_1}{\partial s}(s,0) = 0$, respectively $\frac{\partial F_2}{\partial s}(s,0) = 0$. Moreover, since $F$ is an orientation-preserving  $C^1$-embedding,
 \[ \det(DF(s,x)) > 0 .\]
  At $x = 0$ we have 
\[ \det(DF(s,0)) =\det(  [DF_2(s,0)]_{i, j = 2}^{d-1}), \]
so continuity implies that there is some $\delta >0$ such that 
\begin{equation} \label{det} \det(  [DF_2(s,x)]_{i, j = 2}^{d-1}) > 0 \end{equation}
for every $x \in B^{d-1}(\delta)$.

We define
$$G(s,x)= \left\{ \begin{array}{lcl} F(s,x) & : & s \in [1, 2] \\ (s + F_1(1, x), F_2(1,x)) & : & s \in [3/4,1) \\
(s + 4(s - 1/2)F_1(1, x),  F_2(1, x))& : & s \in [1/2, 3/4)\\
(s  , (4s -1)^{-1}F_2(1, (4s - 1)x))& : & s \in (1/4, 1/2) \\
(1/4, [DF_2(1,0)]_{i,j = 2}^{d} \ x   ) & : & s = 1/4 \\  
(s, [a_{i,j}(s)]_{i,j =1}^{d-1} \ x) & : & s \in [0, 1/4) \\
(s,x) &: & s \in [-1, 0) \\
   \end{array} \right. $$
where $[a_{i,j}(s)]_{i,j = 1}^{d-1}$ is a smooth arc of matrices with positive determinant with 
\[ [a_{i,j}(1/4)]_{i,j = 1}^{d-1} = [DF_2(1,0)]_{i, j = 2}^{d}, \quad [a_{i,j}(0)]_{i,j = 1}^{d-1} = I.\]
(Recall that $\det([DF_2(1,0)]_{i, j = 2}^{d}) > 0$.)  It is easy to check that $G(s,x)$ is a continuous function satisfying (i) -- (iii).

We show $G$ satisfies (iv). On the interval $[1,2]$ injectivity is clear since $F$ is an embedding. 

For the next interval, $[3/4, 1)$ we have that $(1 + F_1(1,x), F_2(1, x))$ is an embedding from $\{ 1 \} \times \R^{d-1} \into \R^d$, by the same reasoning. Suppose that 
$$(s + F_1(1,x), F_2(1, x)) = (s' + F_1(1,x'), F_2(1, x')).$$  
From equation \eqref{det} we have $\det([DF_2(1, x)]_{i,j = 2}^{d}) >   0$, so using the inverse function theorem and possibly decreasing $\delta$, $F_2(1, \cdot)$ is invertible in  $B^{d-1}(\delta)$.  Hence $x = x'$ and then also $s = s'$, so $G$ is injective on $[3/4,1) \times B^{d-1}(\delta)$.  

On $[1/2,3/4)$, as above we have that $F_2(1,x) = F_2(1,x')$ implies that $x = x'$ if $x, x' \in B^{d-1}(\delta)$. Moreover, if $G(s,x) = G(s', x)$, then $ s= s'$, so again $G$ is injective on $[1/2, 3/4) \times B^{d-1}(\delta)$.

In $(1/4, 1/2)$ if $G(s,x) = G(s',x')$ it is automatic that $s = s'$, hence the result follows again since $F_2(1,x)$ is injective on $B^{d-1}(\delta)$.

At $s = 1/4$ the result follows from the Jacobian condition on $F_2$ (equation \eqref{det}).

For the interval $(0, 1/4)$, the result follows from the fact that each $[a_{i,j}(s)]_{i,j=1}^{d-1}$ have determinant greater than zero, hence $G(s,x)$ is injective on $[-1, 2] \times B^{d-1}(\delta)$.
\end{proof}
 \en

We remark that  $G(s,x)$ is $C^1$ except for when $s \in \{1/2, 3/4, 1\}$. With a bit more care, one can show that $G(s,x)$ could be made $C^1$, but we will not need this.

\bn  \begin{nlemma}
\label{Z:5}
There exists $\epsilon > 0$ and 
a continuous injection 
\[
\tau_{0} : [-\epsilon, 1 + \epsilon] \times 
\overline{B^{d-1}(1)} \rightarrow S^{d}
\] 
satisfying the following.
\begin{enumerate}
\item
For $| s | \leq \epsilon$ and $x$ in $\overline{B^{d-1}(1)}$, \[
\tau_{0}(s+1, x) = \varphi(\tau_{0}(s, x)).
\]
\item For all $n \neq 0$, 
\[
\varphi^{n}( \tau_{0} ([0, 1 ) \times \{ 0 \}) ) 
\cap \tau_{0} ([0, 1 ) \times \{ 0 \}) = \emptyset.
\]
\end{enumerate}
\end{nlemma}

\begin{proof} Let $\lambda_{\R} : \R \to S^d$ be as defined following the proof of Lemma~\ref{Z:3}.  By the tubular neighbourhood theorem there exists a $C^1$ embedding 
\[  \alpha : [-1/2, 1 + 1/2] \times \R^{d-1} \to S^d \]
such that 
\begin{equation} \label{alpha} 
\alpha(s, 0) = \lambda_{\R}|_{[-1/2, 1+ 1/2]}.
\end{equation}
There exist $\epsilon, \gamma > 0$ such that 
\[ \varphi \circ \alpha([-\epsilon, \epsilon] \times B^{d-1}(\gamma)) \subset \textup{Im} (\alpha). \]

Let $\beta : \R \times B^{d-1}(\gamma) \to \R \times B^{d-1}(\gamma)$ be given by $\beta(s, x) = (s-1, x)$. Define $F : [1-\epsilon, 1 + \epsilon] \times B^{d-1}(\gamma) \to \R^d$ by
\[ F = \alpha^{-1} \circ \varphi \circ \alpha \circ \beta .\]
By construction, $F$ is an orientation-preserving $C^1$-embedding satisfying $F(s, 0) = (s, 0)$ for every $s \in [1-\epsilon, 1 + \epsilon]$. Hence we may apply the previous lemma to get 
$$G : [-\epsilon, 1+ \epsilon] \times B^{d-1}(\gamma) \to [-1/2, 1 + 1/2] \times \R^{d-1}.$$  (Note that the intervals are slightly different but this does not matter). From the lemma we have that $G(s, 0) = (s, 0)$, for $s \in  [-\epsilon, 1+ \epsilon]$ and that $G$ is injective by possibly shrinking $\gamma$. 

Define $\tau_0 : [-\epsilon, 1 + \epsilon] \times \overline{B^{d-1}(1)} \to S^d$ by $\tau_0 = \alpha \circ G$, (where we are tacitly using the fact that $\overline{B^{d-1}(1)}$ can be identified with  $\overline{B^{d-1}(\gamma/2)}$.)

We now show that $\tau_0$ has properties (i) and (ii). Property (ii) follows immediately from the fact that $\tau_0(s, 0) = \lambda_{\R}(s)$ and Lemma~\ref{Z:3} (iii). For property (i) we have, for $s \in [-\epsilon, \epsilon]$, that 
\[ \tau_0(s+1, x) = (\alpha \circ G) (s + 1, x)  = \varphi \circ \alpha \circ \beta(s+1, x) = \varphi(\alpha(s,x)),\]
and 
\[ \varphi(\tau_0(s, x)) = \varphi (\alpha \circ G (s, x))  = \varphi (\alpha (s, x)). \]
\end{proof} 
\en

Given $\tau_0$ as in the previous lemma, we define $ \tau : \R \times \overline{B^{d-1}(1)} \to S^d $ by
\begin{equation} \label{tau}
\tau(s, x) = \varphi^{\lfloor s \rfloor}\tau_0(s \textup{ mod } 1, x).
\end{equation}

\bn
  \begin{nlemma}
\label{Z:10}
There exists a sequence of positive numbers 
$1 > \rho_{1} > \rho_{2} > \cdots >0$ such that 
$\tau |_{[-n-2, n+2] \times \overline{B^{d-1}(2\rho_{n})}}$
is injective.
\end{nlemma}  
\en

Define functions $r_{n}: \R \rightarrow [0, 1]$ by 
\[
r_{n}(s) = \left\{ \begin{array}{cl}
  0, & | s | > n, \\
  n - | s |, & n- \rho_{n} \leq | s | \leq n, \\
  \rho_{n}, & | s | \leq n - \rho_{n}.
\end{array} \right.
\]
\begin{figure}[t]
\begin{picture}(300,150)
\put(0,20){\line(1,0){300}}
\put(10,0){$-n$}
\put(280,0){$n$}
\put(150,0){$0$}
\put(210,50){$\rho_{n}$}
\put(90,80){$r_{n}(s)$}
\put(200,20){\vector(0,1){50}}
\put(20,20){\line(1,1){50}}
\put(70,70){\line(1,0){160}}
\put(150,20){\line(0,1){100}}
\put(230,70){\line(1,-1){50}}
\end{picture}
\caption{Graph of $r_n(s)$}
\end{figure}
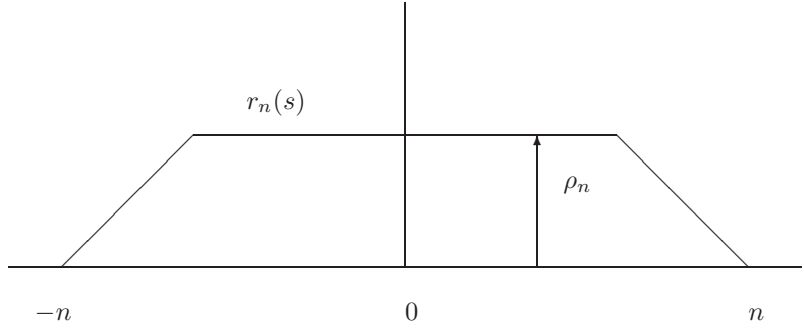

For $n \geq 1$, define
\[
L_{n} = \tau([-n, n] \times \{ 0 \}), \text{ and } L_{\infty} = \tau( \R \times \{ 0 \}).
\]
Also define, for each $n \geq 1$,
\begin{equation} \label{Xn}
X_{n} = S^{d}-  \tau \{ ( s, x) \mid -n \leq s \leq n, \vert x \vert \leq r_{n}(s) \}.
\end{equation}

\bn
  \begin{nlemma}
\label{Z:13}
The closure of $X_{n}$, $\overline{X_{n}}$, is contractible.
\end{nlemma}

\begin{proof}
It is clear that 
$ \tau(\{ ( s, x) \mid -n < s < n, \vert x \vert < r_{n}(s) \})$ is
the complement of $\overline{X_{n}}$, it is an open 
set in $S^{d}$ homeomorphic to an open ball. 
Removing an open set of this form 
in $S^{d}$ yields a space homeomorphic to a closed 
ball in $\R^{d}$, which is contractible.
\end{proof}  
\en

We define a function $\beta_{n}: S^{d} - L_{n} \rightarrow S^{d}$ as follows.
Let 
\begin{equation} \label{Rn}
R_{n} = \tau \{ (s, x) \mid -n \leq s \leq n, 0 < \vert x \vert \leq 2r_{n}(s), x \in \R^{d} \}.
\end{equation}
Observe that  any point in $R_{n}$ may be written 
uniquely as $\tau(s, t x)$, with $ -n \leq s \leq n$, 
$0 < t \leq  2r_{n}(s)$ and $x$ in $S^{d-2}$. 
We then define
\begin{equation} \label{betan}
\beta_{n}(\tau(s, tx)) = \tau\left( s, \left( \frac{t}{2} + r_{n}(s)\right) x \right) .
\end{equation}
Observe that $\beta_{n}$ fixes any 
point with  $t  = 2r_{n}(s)$.
We set $\beta_n$ to be the identity on 
$S^{d} - L_{n} - R_{n}$.

 \bn \begin{nlemma}
\label{Z:15} For $\beta_n : S^d - L_n \to S^d$ defined as in \eqref{betan}, the following hold:
\begin{enumerate}
\item $\beta_{n}$ is continuous.
\item The image  $\beta_{n}( S^{d} - L_{n} )$ is $ X_{n}$.
\item $\beta_{n} $ is injective.
\item $\beta_n^{-1}: X_{n} \rightarrow S^{d} - L_{n}$ is continuous.
\end{enumerate}
\end{nlemma} 
\en

\begin{proof}
Each of (i) -- (iv) is a straightforward calculation. 
\end{proof}

Let us add one more useful observation: at this point, $\beta_{n}$ is 
not defined on $\tau(\pm n, 0)$, but if we extend the definition to leave these
fixed, it is still continuous.

Let $d_{0}$ be any fixed  metric on $S^{d}$ which yields the usual topology.
Of course, since $S^{d}$ is compact, $d_{0}$ is bounded.

We define a sequence of metrics, $d_{n}$ on $S^{d} - L_{n}$, by 
\begin{equation} \label{dn}
d_{n}(x, y) = d_{0}(\beta_{n}(x), \beta_{n}(y)), \quad  x, y \in S^{d} - L_{n}.
\end{equation}
That is, 
$\beta_{n}$ is an isometry from $ (S^{d} - L_{n}, d_{n})$ to $(X_{n}, d_{0})$.

 \bn \begin{ndefn} Define $Z_{n}$ to be the completion of $S^{d} - L_{n} $ in 
$d_{n}$.
\end{ndefn}  
\en

Given these definitions, the following is obvious since the completion of 
$X_{n}$ in the metric $d_{0}$ is simply its  closure.

\bn
  \begin{nlemma}
\label{Z:20}
The map $\beta_{n}$ extends to a homeomorphism 
from $(Z_{n}, d_n)$ to $(\overline{X_{n}}, d_0)$. In particular,
$Z_{n}$ is a connected, contractible, compact absolute neighbourhood retract with 
$\dim(Z_{n}) = d$.
\end{nlemma}  
\en
 
Now we establish the following important relations between our metrics. 
 
 \bn
  \begin{nlemma}
\label{Z:25}
Let $(z_k)_{k \in \mathbb{N}}$ be a sequence in $S^{d} - L_{n}$, $n \geq 1$. Then
\begin{enumerate} 
 \item  if $(z_k)_{k \in \mathbb{N}}$ is Cauchy in $d_{n}$ then it is also Cauchy in $d_{0}$,
\item  if $(z_k)_{k \in \mathbb{N}}$ is Cauchy in $d_{n}$  then it is also Cauchy in $d_{n-1}$.
 \end{enumerate}
\end{nlemma}

\begin{proof}
Let us consider the first part. It suffices to prove that 
if $(z_k)_{k \in \mathbb{N}}$ is any sequence in $S^{d} - L_{n}$ 
such that $\beta_{n}(z_{k})$ is Cauchy in $S^{d}$ in $d_{0}$,
then $(z_k)_{k \in \mathbb{N}}$ itself is Cauchy in $d_{0}$.

 Let us first suppose that $(z_k)_{k \in \mathbb{N}}$
lies entirely in   the complement of $R_{n}$.
Then $\beta_{n}$ is the identity on $z_{k}$ and the conclusion 
is clear. Now let us assume the sequence lies entirely inside
 $R_{n}$.
We write  $z_{k} = \tau((s_{k}, t_{k} x_{k}))$, with
 $s_{k} \in [-n, n]$,
$0 < t_{k} \leq 2 r_{n}(s_{k})$ and $x_{k}$ in $S^{d-2}$.
 The fact $(\beta_{n}(z_{k}))_{k \in \mathbb{N}}$
is  Cauchy in $d_{0}$ and that $\tau$ is uniformly continuous
means that $\left( s_{k}, \left( \frac{t_{k}}{2} + r_{n}(s_{k}) \right) x_{k} \right)$ 
is Cauchy in the usual metric of $[-n, n] \times \R^{d-1}$.
It follows that $(s_{k})_{k \in \N}$ converges to some $s$ in 
$[-n, n]$, while $ \left( \frac{t_{k}}{2} + r_{n}(s_{k}) \right) x_{k}$ 
converges to $y$. Taking norms and recalling that $x_{k}$ is a unit 
vector, we see that $\frac{t_{k}}{2} + r_{n}(s_{k})$ converges
to $\vert y \vert$. Putting these together we have 
\[
2 \vert y \vert = \lim_{k} t_{k} + 2 r_{n}(s_{k}) =  2 r_{n}(s) + \lim_{k} t_{k}
\]
and so $t_{k}$ converges to $2 (\vert y \vert - r_{n}(s)) $.

If  $\vert y \vert - r_{n}(s) = 0$, then $t_{k} x_{k}$ converges to the 
zero vector. If  $\vert y \vert - r_{n}(s) \neq 0$, then 
$x_{k}, k \geq 1$ converges to 
$\left( 2 (\vert y \vert - r_{n}(s)) \right)^{-1}y$ and again
$t_{k} x_{k}$ is convergent. 

Finally, we need to consider the case where $(z_k)_{k \in \mathbb{N}}$ contains
infinitely many terms outside the region and infinitely many terms 
inside the region. From the arguments above, the
 two subsequence each converge 
in the usual topology to two points, say $z_{out}$ and $z_{in}$, respectively.
But we also have $\beta_{n}(z_{out}) = z_{out}$. So the entire sequence
$\beta_{n}(z_{k})$ is converging to the point $z_{out}$ which lies 
in the range of $\beta_{n}$. Then the conclusion follows from 
the fact that $\beta_{n}^{-1}$ is continuous. This completes the proof 
of the first statement.

For the second part,  it suffices to prove that 
if $(z_k)_{k \in \mathbb{N}}$ is any sequence in $S^{d} - L_{n}$ 
such that $(\beta_{n}(z_{k}))_{k \in \N}$ is Cauchy in $S^{d}$ in the usual metric,
then $(\beta_{n-1}(z_{k}))_{k \in \N}$  also Cauchy.

Let us consider the case when the 
sequence $(z_k)_{k \in \mathbb{N}}$ is outside  $R_{n}$. In this case,
we have $\beta_{n}(z_{k}) = z_{k}$.
 Then $(z_k)_{k \in \mathbb{N}}$ is 
Cauchy and hence convergent. We note that the limit point lies 
outside of $R_{n}$ or on its boundary, and hence is not in 
$L_{n-1}$. It follows since $\beta_{n-1}$ is 
continuous that $\beta_{n-1}(z_{k})$ is convergent and hence Cauchy.

Now suppose that $z_{k}$ is in $R_{n}$ and write 
$z_{k} = \tau(s_{k}, t_{k} x_{k})$ as before. As in the first case, we 
know that $s_{n}$ converges to some $s$ in $[-n, n]$, 
 while $ \frac{t_{k}}{2} + r_{n}(s_{k}) x_{k}$ 
converges to $  y $. 
If $\vert y \vert - r_{n}(s) \neq 0$, then as before, 
$z_{k}$ is converging to a point in $S^{d} - L_{n}$,
 which is a subset of $S^{d}-L_{n-1}$. We use the fact that
 $\beta_{n-1}$ is continuous to conclude that $\beta_{n-1}(z_{k})$  
 is convergent and hence Cauchy. 

Now suppose that $\vert y \vert - r_{n}(s) = 0$. In this case, $t_{k}$ is converging 
to zero and this means that 
\[
y = \lim_{k}\left( \frac{t_{k}}{2} + r_{n}(s_{k}) \right) x_{k}
  = \lim_{k} r_{n}(s_{k}) x_{k}.
  \]

First suppose that $ 1-n < s < n-1  $. This implies that 
$r_{n}(s) = r_{n} \neq 0$ and so $x_{k}$ is actually convergent. It follows that  
\[
\lim_{k} \beta_{n-1}(z_{k}) = \lim_{k} \tau\left( s_{k}, \left(\frac{t_{k}}{2} + r_{n-1}(s_{k}) \right) x_{k} \right)
= \tau\left( s, r_{n-1}(s) \lim_{k} x_{k} \right)
\]
 and so $\beta_{n-1}(z_{k}), k \geq 1$ is Cauchy.
 
 Next, suppose that $n-1 \leq |s | \leq n$. Here we have  $r_{n-1}(s) = 0$ and so
 \[
\lim_{k} \beta_{n-1}(z_{k}) = \lim_{k} \tau\left( s_{k}, \left(\frac{t_{k}}{2} + r_{n-1}(s_{k}) \right) x_{k} \right)
= \tau( s, 0)
\]
and again $(\beta_{n-1}(z_{k}))_{k \in \N}$ is Cauchy.      

Finally, suppose $(\beta_{n-1}(z_{k}))_{k \in \N}$ contains infinitely many terms outside the region and infinitely many terms inside the region. As in the proof  of (i), we can find two subsequences $(\beta_{n-1}(z_{k}))_{\textup{in}}$ and $(\beta_{n-1}(z_{k})_{\textup{out}})$. We know $(\beta_{n-1}(z_{k})_{\textup{out}})$ converges to some $x_{\textup{out}}$ by the above. One checks that $x_{\textup{out}}$ lies outside $R_n$ so $\beta_n(x_{\textup{out}}) = x_{\textup{out}}$. By the same reasoning as in part (i), we have that $(\beta_n(z_k))_{k \in \N}$ converges to $x_{\textup{out}}$. Then by continuity of $\beta_n, \beta_n^{-1}$ and $\beta_{n-1}$ we have that 
\[\beta_{n-1}(z_k) = \beta_{n-1} \circ \beta_n^{-1} \circ \beta_n (z_k) \to  \beta_{n-1}(x_{\textup{out}}), \text{ as } k \to \infty \]
so $(\beta_{n-1}(z_k))_{k \in \N}$ is also Cauchy.
\end{proof}  \en

Observe that, for every $n \geq 1$,
 $S^{d}$ is the completion
of $S^{d} - L_{n} $ in $d_{0}$. It is also the completion 
of $S^{d} - L_{\infty}$ in $d_{0}$.

The immediate consequence of the last lemma 
is that the identity map
on $S^{d} - L_{\infty}$ extends to  well-defined,
continuous maps 
$\pi_{n}: Z_{n} \rightarrow Z_{n-1}$ and 
$q_{n}: Z_{n} \rightarrow S^{d}$.  Furthermore, both $\pi_n$ and $q_n$ are surjective since in both cases $\textup{Im}(Z_n)$ is compact and dense.

Define a metric $d_{\infty}$ on $S^{d} - L_{\infty}$ by 
\[
d_{\infty}(x, y) = \sum_{n \geq 1} 2^{-n} d_{n}(x, y).
\]
This uses the fact that each $d_{n}$ is bounded by the same number which 
bounds $d_{0}$. We observe that a sequence 
$(z_{k})_{k \in \N}$ is Cauchy in $d_{\infty}$ 
if and only if it is Cauchy in each $d_{n}$.

\bn  \begin{ndefn}
Define $ Z$ to be the completion of  $S^d - L_{\infty}$ in the metric $d_{\infty}$.
\end{ndefn}  \en

Based on the relationships between the $(Z_n, d_n)$, $n \geq 1$, we also have the following description of $(Z, d_{\infty})$: 

\bn  \begin{nlemma}
\label{Z:30}
The space $Z$ is homeomorphic to the inverse limit of the system 
\[ Z_{1} \stackrel{\pi_{2}}{\leftarrow} Z_{2} \stackrel{\pi_{3}}{\leftarrow} \cdots. \]
\end{nlemma} \en
 
Many of the properties of $Z_{n}, n \geq 1$, which we observed in Lemma \ref{Z:20} are preserved under inverse limits:  compactness, 
connectedness and finite covering dimension. Of course, contractibility is not preserved, but both cohomology and $K$-theory are continuous. Therefore we have the following:

\bn
 \begin{ncor}
\label{Z:35}
The space $Z$ is an infinite connected, compact metric space with finite covering dimension. It has the same cohomology and $K$-theory as a point; in fact, $C(Z)$ is $KK$-equivalent to $\mathbb{C}$.
\end{ncor} 

\begin{proof}
The only thing we have not shown is the $KK$-equivalence, but this follows from the UCT and $K_*(C(Z)) \cong \mathbb{Z} \oplus 0$.
\end{proof}
\en
 
We can even say more about the covering dimension of $Z$:

\bn
\begin{nlemma} \label{dimZ}
The covering dimension of $Z$ is $d$ or $d-1$.
\end{nlemma}

\begin{proof}
Since $\dim(Z_n) = d$ for every $n \in \mathbb{N}$, it is clear that $\dim(Z) \leq d$. To show that it is at least $d-1$, we exhibit an embedding $\sigma : [-1/2, 1,2] \times S^{d-2} \to Z$. It is enough to show that there are embeddings $\sigma_n :  [-1/2, 1,2] \times S^{d-2} \to Z_n$ for every $n$ which are consistent with the inverse limit. Define
\[ \sigma_n : [-1/2, 1,2] \times S^{d-2} \to Z_n , \quad (s,x) \mapsto \lim_{k \to \infty} \tau(s, k^{-1}x) .\]

To show that $\sigma_n $ is well defined, we show that $\tau((s, k^{-1}x)_{k\in\mathbb{N}})$ is Cauchy with respect to the metric $d_n$.

Take $K \in \mathbb{N}$ large enough so that $\tau((s, k^{-1}x)) \in R_n$ for every $k \geq K$. Then we have 
\begin{eqnarray*}
d_n( \tau((s, k_1^{-1}x)), \tau((s, k_2^{-1}x)) &=& d_0(\beta_n( \tau((s, k_1^{-1}x))), \beta_n(\tau((s, k_2^{-1}x))) \\
&=& d( (k_1^{-1}/2 + \rho_n)x, (k_2^{-1}/2 + \rho_n)x)\\
&=& 1/2 | k_1^{-1} - k_2^{-1} |,
\end{eqnarray*}
from which we see that the sequence is Cauchy with respect to $d_n$.

Now we must show that $\sigma_n$ is injective. Suppose that $(s,x) \neq (s',x')$ and let $K$ be sufficiently large so that $\tau((s,x)), \tau((s',x')) \in R_n$. Suppose that $s \neq s'$. Then for every $k \geq K$,
\begin{eqnarray*}
d_n( \tau((s, k_1^{-1}x)), \tau((s', k^{-1}x')) &=& d_0(s, (k^{-1}/2 + \rho_n)x, (s', k^{-1}/2 + \rho_n)x')\\
&\geq& | s - s' | \\
&>& 0. 
\end{eqnarray*}

On the other hand, if $ s = s'$ we have that $x \neq x'$, and in this case, for all $k \geq K$, we have
\begin{eqnarray*}
d_n( \tau((s, k_1^{-1}x)), \tau((s, k^{-1}x')) &=& d_0(\beta_n(s, (k^{-1}/2 + \rho_n)x), \beta_n((s, k^{-1}/2 + \rho_n)x'))\\
& \geq & d_{S^{d-2}}(\rho_n x, \rho_n x') \\
&>& 0,
\end{eqnarray*}
where $d_{S^{d-2}}$ denotes the metric on $S^{d-2}$.

Finally, since for every $k \in \mathbb{N}$, we have that $\sigma_n((s, x)) \in S^d - L_{\infty}$ so it is clear that $\pi_n \circ \sigma_n = \sigma_{n-1} \circ \pi_n$.
\end{proof}
\en

The identity map
on $S^{d} - L_{\infty}$ extends to a  well-defined, 
continuous, surjective  map
$q: Z \rightarrow S^{d}$. 

\bn
  \begin{nlemma}
\label{Z:28}
For each point $x$ in $S^{d} - L_{\infty}$, $q^{-1}\{ x \}$
is a single point.
\end{nlemma}

\begin{proof}
We treat the points of $Z$ and $S^{d}$ as equivalence classes of 
Cauchy sequences in $S^{d} - L_{\infty}$ 
in the metrics $d_{\infty}$ and $d_{0}$, 
respectively. Let $(w_{k})_{k \in \N}$ and $(z_{k})_{k \in \N}$ 
be two Cauchy sequences $S^{d} - L_{\infty}$  in $d_{\infty}$ 
and suppose they map to the same point under $q$, That is, they 
are equivalent in the metric $d_{0}$. Suppose also 
that they converge to a point $z$ in $S^{d} - L_{\infty}$ 
in $d_{0}$. Fix $n \geq 1$. As $z$ is in 
$S^{d} - L_{\infty}$ it is also in $S^{d} - L_{n}$. The latter is open in $S^{d}$ in the metric $d_{0}$, so we may find an 
open ball $B$ whose closure is contained in 
 $S^{d} - L_{n}$. For $k$ sufficiently large, both 
 $z_{k}$ and $w_{k}$ are in $B$. The function 
 $\beta_{n}$ is defined and  
 continuous on $B$ and so we conclude that 
 both sequences $\beta_{n}(z_{k})$ and $\beta_{n}(w_{k})$ 
 are converging to $\beta_{n}(z)$ in $d_{0}$. 
 Hence, the sequences $z_{k}$ and $w_{k}$ are 
 equivalent in the $d_{n}$ metric. As this is true for every 
 $n$, these sequences are also equivalent in the 
 $d_{\infty}$ metric and we are done. 
\end{proof}  
\en

We now turn to the problem of defining our minimal dynamical system.

\bn
  \begin{nlemma}
\label{Z:40}
If a sequence $(z_k)_{k \in \mathbb{N}}$ in $S^{d} - L_{\infty} $ is Cauchy in 
$d_{n}$, with $n \geq 2$,
 then both sequences  $(\varphi(z_{k}))_{k \in \N}$ and 
 $(\varphi^{-1}(z_{k}))_{k \in \N}$ are Cauchy in $d_{n-1}$.
\end{nlemma}

\begin{proof}
We consider only the case $(\varphi(z_{k}))_{k \in \N}$, the other being similar.
Once again it suffices to assume that $(\beta_{n}(z_{k}))_{k \in \N}$ is Cauchy
in $d_{0}$  and show that $(\beta_{n-1}\circ \varphi(z_{k}))_{k \in \N}$ is also Cauchy in $d_0$.

If the entire sequence lies outside of 
$R_{n}$, then $z_{k} = \beta_{n}(z_{k})$ converges (in the $d_0$ metric) to some point, say $y$,
in the complement of $R_{n}$ or on its boundary. If $y$ is not in $L_{n}$, then 
$\varphi(y)$ is not in $L_{n-1}$ and then same argument using the 
continuity of $\beta_{n-1}$ gives the desired result. 
The only point  in the closure
of the complement of $R_{n}$ with $\varphi(y)$ in $L_{n-1}$ is $\tau(-n, 0)$
with 
$\varphi(y) = \tau(1-n,0)$. We noted earlier that $\beta_{n-1}$ extends continuously to
this point by fixing it and so the same argument works here.

Now assume that $z_{k}$ lies in $R_{n}$ but the sequence $\varphi(z_{k})$ lies outside $R_{n-1}$. In this case,  $\beta_{n-1}(\varphi(z_k)) = \varphi(z_k)$.
But we know from the first part of Lemma~\ref{Z:25} that $(z_{k})_{k\in \mathbb{N}}$ is Cauchy with respect to $d_0$
and hence so is $\varphi(z_{k})$ and we are done.

Now assume $z_{k}$ lies in $R_{n}$ and $ \varphi(z_{k})$ lies in $R_{n-1}$. As before we write 
\[ z_k = \tau(s, tx) ,\]
with $s \in [-n, n]$, $0< t \leq 2r_n(s)$ and $x \in S^{d-2}$.  
In this case, using the definition of $\tau$ given in \eqref{tau}, we have
\[
 \beta_{n-1} \circ \varphi(z_{k}) =  \tau\left(s_{k} + 1, \left( \frac{t_{k}}{2} + r_{n-1}(s_{k} +1)
\right) x_{k} \right).
\] 
The argument proceeds exactly as before. We know $s_{k}$ converges to $s$
and $t_{k}$ converges to $2( \vert y \vert - r_{n})$.  If this is positive then  $x_{k}$ also converges and the desired result follows easily
from the formula above for $\beta_{n-1}\circ \varphi(z_{k})$. If $\vert y \vert = r_{n}(s)$ then $t_{k}$ converges to zero. We break this up into three cases.

\underline{Case 1}: $s > n-2$. Then $\varphi(z_{k})$ is not in $R_{n-1}$ and this case is already done.

\underline{Case 2}: $-n < s \leq n-2$. Here $r_{n}(s) > 0$
and thus $\left( \frac{t_{k}}{2} + r_{n}(s_{k}) \right) x_{k}$ converging to $y$ implies that the sequence $x_{k}$ itself is convergent.
The convergence of $\beta_{n-1} \circ \varphi(z_{k})$ follows
from the formula above.

\underline{Case 3}:  $s = -n$. In this case, both $t_{k}$ and $r_{n-1}(s_{k}+1)$ 
are converging to  $0$ and the convergence of $\beta_{n-1} \circ \varphi(z_{k})$ 
again follows from the formula above.

Finally, when $(z_k)_{k \in \mathbb{N}}$ has infinitely many terms lying both outside and inside the region, the proof is similar to previous calculations.
\end{proof}  
\en

\bn
  \begin{ncor}
\label{Z:45}
The map $\varphi$ on  $S^{d} - L_{\infty}$ extends to a homeomorphism 
of $(Z, d_{\infty})$,  denoted by $\zeta$. Moreover, we have 
$q \circ \zeta = \varphi \circ q$; that is $q$ is a factor map
from $(Z, \zeta)$ to $(S^{d}, \varphi$).
\end{ncor}  
\en

\bn
  \begin{ntheorem}
\label{50}
The homeomorphism $\zeta$ of $Z$ is minimal. 
\end{ntheorem}

\begin{proof}
Let $Y$ be a non-empty, closed (hence compact), $\zeta$-invariant subset of $Z$. It follows that $q(Y)$ is a non-empty, compact (hence closed) $\varphi$-invariant subset of $S^{d}$ and hence $q(Y) = S^{d}$. As the quotient map $q$ is injective on $S^{d} - L_{\infty} \subseteq Z$, it follows that $S^{d} - L_{\infty} \subseteq Y$. As $Y$ is closed and $Z$ 
is defined as the completion of $S^{d} - L_{\infty} $, we conclude that $Y = Z$ and so $\zeta$ is minimal.
\end{proof}
 \en

 \bn
\begin{ntheorem} \label{measures}
 The factor map $q$ defined in Corollary~\ref{Z:45} induces an affine bijection between the $\zeta$-invariant Borel probability measures on $Z$ and the $\varphi$-invariant Borel probability measures on $S^d$.
\end{ntheorem}

\begin{proof}
Let $\mu$ be  a  $\zeta$-invariant measure.
Then $q^{*}(\mu)$ is $\varphi$-invariant. 
The set $L_{\infty}$ is 
a Borel subset of $S^{d}$. It is also $\varphi$-invariant.
Moreover, the system $\varphi$, restricted to 
$L_{\infty}$ is conjugate to the map $x \rightarrow x+1$ 
on $\R$, which has no nonzero finite invariant measures.
This implies that $q^{*}(\mu)(L_{\infty}) = 0$. This means that 
$S^{d} - L_{\infty}$ has full measure under $\mu$ and as
$q$ is a bijection on this set, the conclusion follows.
\end{proof}  
\en
 
We remark that  one may modify the construction of the embedding 
of $\mathbb{R}$ into $S^d$ to an embedding of two disjoint 
copies of $\mathbb{R}$. Proceeding  in an analogous fashion, 
the space $Z_n$ is homeomorphic to the sphere with two 
open balls removed. This can easily seen to be 
homeomorphic to $[0,1] \times S^{d-1}$. Continuing, 
the space $Z$ can be seen to have the same cohomology
as the even sphere $S^{d-1}$ and admits 
a minimal
homeomorphism (which can be arranged to be uniquely ergodic) while $S^{d-1}$ does not by the Lefschetz--Hopf Theorem.

Using $\sigma$ from the proof of 1.13, it would seem that the map sending $x$ in $S^{d-2}$ to $\sigma(0, x)$ in $Z$ is not homotopic to a constant function and hence $\pi_{n-2}(Z)$ is non-trivial. If this is correct, then $Z$ is not contractable.

\section{A dynamical presentation of the Jiang--Su algebra} \label{Jiang-Su}

We begin this section by recalling some facts about the Jiang--Su algebra (see \cite{JiaSu:Z}).  In what follows, for any $n \in \mathbb{N}$, we let $M_n$ denote the $\mathrm{C}^*$-algebra of $n \times n$ matrices over $\mathbb{C}$. 

 \bn
\begin{ndefn}
Let $p, q \in \mathbb{N}$. 
The $(p,q)$-dimension drop algebra $A_{p,q}$ is the defined to be
\[ A_{p,q} = \{ f \in C([0,1], M_p \otimes M_q) \mid f(0) \in M_p \otimes 1_q, f(1) \in 1_p \otimes M_q  \}.\]
\end{ndefn}
\en
 
Note that when $p$ and $q$ are relatively prime, $A_{p,q}$ is 
projectionless, that is, its only projections are $0$ and 
$1$.

\bn  
\begin{ntheorem}\cite[Theorem 4.5]{JiaSu:Z} Let $G$ be an inductive limit of a sequence of finite cyclic groups and $\Omega$ a nonempty metrizable Choquet simplex. Then there exists a simple unital infinite-dimensional projectionless $\mathrm{C}^*$-algebra $A$ which is isomorphic to an inductive limit of dimension drop algebras and satisfies 
\[ ((K_0(A), K_0(A)_+, [1_A]), K_1(A), T(A)) \cong ((\mathbb{Z}, \mathbb{Z}_+, 1), G, \Omega). \]
\end{ntheorem}
 \en

In the same paper, Jiang and Su showed that any two such simple inductive limits of finite direct sums of dimension drop algebras are isomorphic if and only if their Elliott 
invariants are isomorphic \cite[Theorem 6.2]{JiaSu:Z}. Moreover, the isomorphism of $\mathrm{C}^*$-algebras can be chosen to induce the isomorphism at the level of the invariant. 

\bn 
\begin{ndefn} The Jiang--Su algebra $\mathcal{Z}$ is the unique simple unital infinite-dimensional inductive limit of finite direct sums of dimension drop algebras satisfying 
\begin{eqnarray*}
((K_0(\mathcal{Z}), K_0(\mathcal{Z})_+, [1_{\mathcal{Z}}]), K_1(\mathcal{Z}), T(\mathcal{Z})) &\cong &((\mathbb{Z}, \mathbb{Z}_+, 1), 0, \{ \textup{pt} \} )\\
 & \cong & ((K_0(\mathbb{C}), K_0(\mathbb{C})_+, [1_{\mathbb{C}}]), K_1(\mathbb{C}), T(\mathbb{C}))  
 \end{eqnarray*}
 \end{ndefn}
 \en

The goal of this section is to  exhibit the Jiang--Su algebra $\mathcal{Z}$ as the $\mathrm{C}^*$-algebra of a minimal 
\'{e}tale  equivalence relation.  As described in the introduction, this should be seen in analogy to the von Neumann algebra--measurable dynamical setting where the hyperfinite $\mathrm{II}_1$ factor $\mathcal{R}$ is shown to be the von Neumann algebra an amenable measurable equivalence relation.

 \bn
\begin{ndefn} Let $X$ be a compact metrizable space. Then an equivalence relation $\mathcal{E} \subset X \times X$ with countable equivalence classes is called minimal if every equivalence class is dense in $X$.
\end{ndefn}
\en

Let $(X, \alpha)$ be a minimal dynamical system of an infinite compact metric space. Let $\mathcal{E}_{\alpha} \subset X \times X$ denote the orbit equivalence relation of $(X, \alpha)$.  
As described in \cite{Ren:groupoid}, it is equipped with a natural topology in which it is \'{e}tale. Note that the orbit equivalence relation from a minimal dynamical system is a minimal equivalence relation.

 \bn
\begin{ndefn} For $y \in X$ the orbit-breaking equivalence relation $\mathcal{E}_{y}$ is defined as follows: If $(x,x') \in \mathcal{E}_{\alpha}$ then $(x,x')  \in \mathcal{E}_{y}$ if $\alpha^n(x) \neq y$ for any $n \in \mathbb{Z}$ or there are $n, m > 0$ such that $\alpha^n(x) = \alpha^m(x') = y$ or there are $n, m \leq 0$ such that $\alpha^n(x) = \alpha^m(x') = y$. 
\end{ndefn}
\en

Note that this splits any equivalence class in $\mathcal{E}_{\alpha}$ containing the point $y$ into two equivalence classes: one consisting of the forward orbit, the other of the backwards orbit. It is easily seen to be an open 
subset of $\mathcal{E}$ in 
the relative topology and, with that topology, is also 
\'{e}tale.
 
 \bn
\begin{nremarks}
(i) The $\mathrm{C}^*$-algebra $\mathrm{C}^*(\mathcal{E}_{y})$ is a well-known $\mathrm{C}^*$-subalgebra of the crossed product $C(X) \rtimes_{\alpha} \mathbb{Z}$ generated by  $\{ C(X), uC_0(X - \{y\}) \}$, where $u \in C(X) \rtimes_{\alpha} \mathbb{Z}$ denotes the unitary satisfying $ufu^* = f \circ \alpha^{-1}$ for every $f \in C(X)$. Such subalgebras were originally introduced in the study of the $\mathrm{C}^*$-algebras associated to Cantor minimal systems \cite{Putnam:MinHomCantor}. 

(ii) By abuse of notation, let $\alpha$ also denote the action of $\mathbb{Z}$ on $C(X)$. Comparing definitions, it is easy to see that  $\mathrm{C}^*(\mathcal{E}_{y})$ can also be realised as the covariance algebra of a partial action on $C(X)$ given by the triple $(\alpha, C_0(X - \{y\}), C_0(X - \{\alpha(y)\})$, in the sense of Exel \cite[Definition 3.1, Definition 3.7]{Exel:ParAct}. Thus Theorem~\ref{ZTheorem} below also gives a characterisation of $\mathcal{Z}$ as a $\mathrm{C}^*$-algebra coming from a partial action of $\mathbb{Z}$.
\end{nremarks}
\en

This  next result is well-known, but we rephrase it in terms of equivalence relations and give a proof for completeness.

 \bn
\begin{nprop} \label{simple}
Let $(X, \alpha)$ be a minimal dynamical system on an infinite compact metrizable space. For any $y \in X$, $\mathcal{E}_{y}$ is minimal and $\mathrm{C}^*(\mathcal{E}_{y})$ is simple.
\end{nprop}

\begin{proof}
Since $\alpha$ is minimal, for any point $x \in X$ both the forward orbit and backwards orbit are dense in $X$. It follows that every equivalence classes of $\mathcal{E}_y$ is dense in $X$. Since $\mathcal{E}_y$ is minimal the associated $\mathrm{C}^*$-algebra $\mathrm{C}^*(\mathcal{E}_{y})$ is simple. (That $\mathrm{C}^*(\mathcal{E}_{y})$ is simple also shown in \cite[Proposition 2.5]{LinPhi:MinHom}.)
\end{proof}
 \en

  In what follows, $\Q$ denotes the universal UHF algebra, that is, the UHF algebra with $K_0(\Q) = \mathbb{Q}$.

\bn
\begin{nprop} \label{Zprop}
Let $Z$ be an infinite compact metrizable space with finite covering dimension satisfying $K^0(Z) \cong \mathbb{Z}$ and $K^1(Z) = 0$. Let $\zeta : Z \to Z$ be a  minimal, uniquely ergodic homeomorphism. Then for any $z \in Z$ we have
\[ C^*(\mathcal{E}_{z}) \cong \mathcal{Z} .\]
\end{nprop}

\begin{proof}
First, we claim that the class of the trivial line bundle 
is the generator of $K^{0}(Z)$. Taking any map 
from the one-point space into $Z$ and composing with 
the only map from $Z$ onto a point, the composition
 (in that order)
is clearly the identity. It then follows from our 
hypothesis on $K^{*}(Z)$ that these two maps
actually induce isomorphisms at the level of $K$-theory
and the claim follows. 

Then, by the Pimsner--Voiculescu exact sequence, one calculates that 
\[ K_0(C(Z) \rtimes_{\zeta} \mathbb{Z}) \cong \mathbb{Z} \cong K_1(C(Z) \rtimes_{\zeta} \mathbb{Z}) .\]
Next, we use the six-term exact sequence in \cite[Theorem 2.4]{Put:K-theoryGroupoids} (see also \cite[Example 2.6]{Put:K-theoryGroupoids}) to calculate that
\[ K_0(\mathrm{C}^*(\mathcal{E}_z)) \cong \mathbb{Z}, \qquad  K_1(\mathrm{C}^*(\mathcal{E}_z))  = 0.\]
Furthermore, we have that $T(C(Z) \rtimes_{\zeta} \mathbb{Z}) \cong T(\mathrm{C}^*(\mathcal{E}_z))$ \cite[Theorem
1.2]{LinQPhil:KthoeryMinHoms}. Thus $\mathrm{C}^*(\mathcal{E}_z)$ has the same invariant as $\mathcal{Z}$. We can write $\mathrm{C}^*(\mathcal{E}_z)$ as an inductive limit of $\mathrm{C}^*$-algebras of the form $\mathrm{C}^*(C(Z), C_0(Z - Y_n))$ where the $Y_n$ are nested closed subsets of $Z$ with nonempty interior such that $\cap_{n = 0}^{\infty} Y_n = \{ z \}$ \cite[Remark 2.2]{LinPhi:MinHom}. By \cite[Section 3]{LinQPhil:KthoeryMinHoms}, each $\mathrm{C}^*(C(Z), C_0(Z - Y_n))$ is a subhomogeneous algebra whose base spaces have dimension no more than $\dim(Z)$. Thus $\mathrm{C}^*(\mathcal{E}_z)$ is a simple approximately subhomogeneous $\mathrm{C}^*$-algebra with no dimension growth and, applying \cite[(3.2)]{KirWinter:dr}, finite decomposition rank. 

Since there is only one tracial state, projections separate traces and it follows from \cite[Theorem 1.4]{BlaKumRor:apprdiv} that $\mathrm{C}^*(\mathcal{E}_z) \otimes \Q$ has real rank zero whence $\mathrm{C}^*(\mathcal{E}_z) \otimes \Q$ is tracially approximately finite \cite[Theorem 2.1]{Winter:lfdr}. Now $\mathcal{Z} \otimes \Q \cong \Q$ is AF \cite[Theorem 5]{JiaSu:Z}, hence also TAF. Since both these $\mathrm{C}^*$-algebras are in the UCT class, we may apply \cite[Theorem 5.4]{LinNiu:KKlifting} with \cite[Theorem 3.6]{LinNiu:RationallyTAI} to get that $\mathrm{C}^*(\mathcal{E}_z) \otimes \mathcal{Z} \cong \mathcal{Z} \otimes \mathcal{Z}.$ 

Furthermore,  $\mathrm{C}^*(\mathcal{E}_z) \cong \mathrm{C}^*(\mathcal{E}_z)  \otimes \mathcal{Z}$ since any simple separable unital nonelementary $\mathrm{C}^*$-algebra with finite decomposition rank is $\mathcal{Z}$-stable \cite[Theorem 5.1]{Winter:dr-Z-stable}. Then since $\mathcal{Z} \otimes \mathcal{Z} \cong \mathcal{Z}$ \cite[Theorem 4]{JiaSu:Z}, we conclude that $\mathrm{C}^*(\mathcal{E}_z) \cong \mathcal{Z}.$
\end{proof}
 \en

 \bn
\begin{ntheorem} \label{ZTheorem}
There is a compact metric space $Z$ with minimal,
\'{e}tale  equivalence relation
 $\mathcal{E} \subset Z \times Z$ such that $ \mathrm{C}^*(\mathcal{E}) \cong \mathcal{Z} .$
\end{ntheorem}

\begin{proof}
For any $d >1$ odd, there is a uniquely ergodic diffeomorphism $\varphi : S^d \to S^d$. Following the construction Section~\ref{construct}, there is a minimal dynamical system $(Z, \zeta)$ where $Z$ satisfies the hypotheses of Theorem~\ref{Zprop} by \ref{Z:35}. Hence for any $z \in Z$ we have $\mathrm{C}^*(\mathcal{E}_z) \cong \mathcal{Z}$.
\end{proof}
\en
 
 \bn
\begin{nremark}
For any equivalence relation $\mathcal{E}$ on a compact Hausdorff space $X$, the $\mathrm{C}^*$-algebra $C(X)$ is a Cartan subalgebra of $\mathrm{C}^*(\mathcal{E})$ \cite[Chapter II, Section 4]{Ren:groupoid}. Thus it follows from Lemma~\ref{dimZ} that the Jiang--Su algebra contains an infinite number of nonconjugate Cartan subalgebras.
\end{nremark}
\en

\section{Classification in the non-uniquely ergodic case} \label{classification}

At present, few examples of crossed product $\mathrm{C}^*$-algebras associated to minimal dynamical systems without real rank zero are known. This is largely due to a lack of examples of minimal dynamical systems $(X, \alpha)$ with $\dim X > 0$ and more than one ergodic measure. In general, classification for $\mathrm{C}^*$-algebras without real rank zero is much more difficult. Real rank zero implies a plentiful supply of projections. Not only does this suggest more information is available in the invariant (in particular, the $K_0$-group), but it also makes the $\mathrm{C}^*$-algebras easier to manipulate into a particular form, for example, to show it is tracially approximately finite (TAF) as defined in \cite{Lin:TAF1}.  For a long time, the minimal diffeomorphisms of odd dimensional spheres were the main example of minimal dynamical systems leading to $\mathrm{C}^*$-algebras which were not at least rationally TAF (that is, tracially approximately finite after tensoring with the universal UHF algebra). Their classification remained elusive for quite some time. By Theorem~\ref{50}, our construction gives further examples lying outside the real rank zero case and we are able to use the classification techniques from the setting of the spheres to classify these crossed products. Furthermore, by using Winter's classification by embedding result  \cite[Theorem 4.2]{Win:ClassCrossProd}, we also classify the projectionless $\mathrm{C}^*$-algebras obtained from the corresponding orbit-breaking sub equivalence relations.

For a given minimal homeomorphism $\varphi : S^d \to S^d$ let $Z_{\varphi}$ denote the space constructed in Section~\ref{construct}, and denote by $\zeta$ the resulting minimal homeomorphism $\zeta: Z_{\varphi} \to Z_{\varphi}$ as in Theorem~\ref{50}.

 \bn
\begin{nprop} \label{trace space}
$T(C(Z_{\varphi}) \rtimes_{\zeta} \mathbb{Z}) \cong T(C(S^d) \rtimes_{\varphi} \mathbb{Z}).$
\end{nprop}

\begin{proof}
This follows immediately from Theorem~\ref{measures}, since tracial states on the crossed product are in one-to-one correspondence with invariant Borel probability measures of the dynamical system. 
\end{proof}
 \en

 \bn
\begin{ntheorem} \label{TAICross} As above, $\Q$ denotes the universal UHF algebra. Let $\mathcal{A}$ be the class of simple separable unital nuclear $\mathrm{C}^*$-algebras given by
\[ \mathcal{A} = \{ C(Z_{\varphi}) \rtimes_{\zeta} \mathbb{Z} \mid \varphi : S^d \to S^d,  d> 1 \text{ odd,  is a minimal diffeomorphism  }  \}. \]
Then for any $A \in \mathcal{A}$, $A \otimes \Q$ is tracially approximately an interval algebra (TAI).
\end{ntheorem}

\begin{proof}
The space $Z_{\varphi}$ satisfies the following: it is infinite, compact, connected, finite-dimensional and $U(C(Z_{\varphi})) = U_0(C(Z_{\varphi}))$.  Since $\zeta$ and $\id_{C(Z_{\varphi})}$ induce the same map on $K$-theory of $C(Z_{\varphi})$ it follows from the UCT that $[(\cdot) \circ \zeta^{-1}]  = [\id_{C(Z_{\varphi})}]$ in $KK(C(Z_{\varphi}), C(Z_{\varphi}))$ and hence in $KL(C(Z_{\varphi}), C(Z_{\varphi}))$ \cite[2.4.8]{Ror:encyc}.  It follows from \cite[Theorem 6.1]{Lin:Spheres} that $A \otimes \Q$ is TAI. 
\end{proof}
\en 

 \bn
\begin{ncor} \label{CrossClass}
If $A, B \in \mathcal{A}$ then $A \cong B$ if and only if $T(A) \cong T(B)$.
\end{ncor}

\begin{proof}
The previous theorem and \cite[Corollary 11.9]{Lin:asu-class} imply that $A \otimes \mathcal{Z} \cong B \otimes \mathcal{Z}$ if and only if $\Ell(A \otimes \mathcal{Z}) \cong \Ell(B \otimes \mathcal{Z})$. By \cite[Theorem B]{TomsWinter:PNAS} (or \cite[Theorem 0.2]{TomsWinter:minhom}) $A$ and $B$ are both $\mathcal{Z}$-stable. For any $\varphi: S^d \to S^d$ minimal homeomorphism we have
\[K_0(C(Z_{\varphi}) \rtimes_{\zeta} \mathbb{Z}) \cong \mathbb{Z}, \quad K_1(C(Z_{\varphi} )\rtimes_{\zeta} \mathbb{Z}) \cong \mathbb{Z}, \]
which follows immediately from the Pimsner--Voiculescu exact sequence and the fact that $\zeta_*$ is the identity on $K$-theory. Thus Elliott invariants are the same up to the tracial state space and we see $A \cong B$ if and only if $T(A) \cong T(B)$.
\end{proof}
 \en

Let $z \in Z_{\varphi}$ and let $\mathcal{E}_z$ denote the equivalence relation given by breaking the orbit of $\zeta$ at the point $z$.

 \bn
\begin{ntheorem} For every $z \in Z$, the $\mathrm{C}^*$-algebra $\mathrm{C}^*(\mathcal{E}_z) \otimes \mathcal{Q}$ is TAI.
\end{ntheorem}

\begin{proof}
By \cite[Section 3]{LinQPhil:KthoeryMinHoms} $\mathrm{C}^*(\mathcal{E}_z)$ is an inductive limit of recursive subhomogeneous $\mathrm{C}^*$-algebras with base spaces of dimension less than or equal to $\dim(Z_{\varphi})$. Since $\dim(Z_{\varphi}) < \infty$ it follows that $\mathrm{C}^*(\mathcal{E}_z)$ has finite nuclear dimension. By Proposition~\ref{simple}, $\mathrm{C}^*(\mathcal{E}_z)$ is simple. 

Let $ \iota : \mathrm{C}^*(\mathcal{E}_z) \otimes \Q \into A \otimes \Q $ denote the unital embedding. Then the map induced on the tracial state spaces $T(\iota) : T(\mathrm{C}^*(\mathcal{E}_z) \otimes \Q) \to T((C(Z_{\varphi}) \rtimes_{\zeta} \mathbb{Z} ) \otimes \Q)$
is a homeomorphism by Proposition~\ref{trace space} and the fact that $\Q$ has a unique tracial state $\tau_{\Q}$. Moreover,
\[ \iota_0 : K_0(\mathrm{C}^*(\mathcal{E}_z) \otimes \Q) \to K_0(C(Z_{\varphi}) \rtimes_{\zeta} \mathbb{Z} ) \otimes \Q) \]
is an ordered group isomorphism \cite[Lemma 4.3]{StrWin:Z-stab_min_dyn}, (see also \cite[Theorem 4.1 (5)]{Phillips:recsub}). 
Since $K_0(C(Z_{\varphi}) \rtimes_{\zeta} \mathbb{Z} )) \cong \mathbb{Z}$, we have that $S(K_0(C(Z_{\varphi}) \rtimes_{\zeta} \mathbb{Z} )))$ is a point.  Thus, for any $\tau, \tau' \in T(C(Z_{\varphi}) \rtimes_{\zeta} \mathbb{Z} ))$, $\tau_* = \tau'_* \in S(K_0(C(Z_{\varphi}) \rtimes_{\zeta} \mathbb{Z} )))$. Moreover, since any tracial state on $C(Z_{\varphi}) \rtimes_{\zeta} \mathbb{Z} ) \otimes \Q$ is of the form $\tau \otimes \tau_{\Q}$, it follows that $\tau_* = \tau'_* \in S(K_0(C(Z_{\varphi}) \rtimes_{\zeta} \mathbb{Z} ) \otimes \Q))$ for any $\tau, \tau' \in T(C(Z_{\varphi}) \rtimes_{\zeta} \mathbb{Z} ) \otimes \Q)$.

By Theorem~\ref{TAICross}, $C(Z_{\varphi}) \rtimes_{\zeta} \mathbb{Z} ) \otimes \Q$ is TAI. It now follows from \cite[Theorem 4.2]{Win:ClassCrossProd} that $\mathrm{C}^*(\mathcal{E}_z)  \otimes \mathcal{Q} \otimes \Q \cong \mathrm{C}^*(\mathcal{E}_z)  \otimes \mathcal{Q} $ is TAI.
\end{proof}
 \en

 \bn
\begin{ncor} \label{DimDrop}
Let 
\[ \mathcal{B} = \{ \mathrm{C}^*(\mathcal{E}_z)   \mid z \in Z_{\varphi},  \varphi : S^d \to S^d, d > 1 \text{ odd, is a minimal diffeomorphism}  \}.\]
 Then $A, B \in \mathcal{B}$ are isomorphic to projectionless inductive limits of prime dimension drop algebras, and $A \cong B$ if and only if $T(A) \cong T(B)$.
\end{ncor}

\begin{proof} If $A \in \mathcal{B}$ then $A$ is $\mathcal{Z}$-stable by \cite{Win:Z-stabNucDim}. After noting this, the proof that $A \cong B$ if and only if $T(A) \cong T(B)$ is as in Corollary~\ref{CrossClass}. For any $z \in Z_{\varphi}$ we have 
\[K_0(\mathrm{C}^*(\mathcal{E}_z)) \cong \mathbb{Z}, \qquad K_1(\mathrm{C}^*(\mathcal{E}_z)) = 0, \qquad T(\mathrm{C}^*(\mathcal{E}_z)) \cong T(C(S^d) \rtimes_{\varphi} \mathbb{Z}),\]
which follows from \cite[Example 2.6]{Put:K-theoryGroupoids}, \cite[Theorem 4.1]{Phillips:recsub} and Proposition~\ref{trace space}. Now it follows from \cite[Theorem 4.5]{JiaSu:Z} that  $A, B \in \mathcal{B}$ are isomorphic to projectionless inductive limits of prime dimension drop algebras.
\end{proof}
 \en

\section{Outlook} \label{outlook}

Although our construction shows that $\mathcal{Z}$ can be realized as a minimal \'etale equivalence relation, it is certainly not unique (we can start with any odd dimensional sphere, for example).  It would be interesting to further investigate the possibility of realizing various properties of $\mathcal{Z}$ at the dynamical level. For example, could there be a suitable notion of ``strongly self-absorbing'' at the level of equivalence relations? Could we see regularity properties such as mean dimension zero (which may be equivalent to $\mathcal{Z}$-stability of the crossed product $\mathrm{C}^*$-algebra) after appropriately taking a product with our system?

\noindent \textbf{Acknowledgements.}  The authors thank the Banff International Research Station and the organizers 
of the workshop ``Dynamics and $\mathrm{C}^*$-Algebras: Amenability and Soficity'', where this project was initiated. In particular, this work followed  a suggestion of Wilhelm Winter, for  which the authors are indebted.  The authors thank the Department of Mathematics and Statistics at the  University of Victoria and the Banach Center, Institute of Mathematics of the Polish Academy of Sciences for funding research visits facilitating this collaboration. The third author thanks Stuart White for useful discussions. Further thanks go to Magnus Goffeng and Adam Skalski for proofreading an initial draft. Finally, thanks to Alcides Buss, Xin Li and Aidan Sims for pointing out Remarks 2.6(ii) and 2.10 following a talk on this work at the ``Workshop on $\mathrm{C}^*$-algebras: Groups and Actions'' at the University of M\"unster.

\providecommand{\bysame}{\leavevmode\hbox to3em{\hrulefill}\thinspace}
\providecommand{\MR}{\relax\ifhmode\unskip\space\fi MR }
\providecommand{\MRhref}[2]{%
  \href{http://www.ams.org/mathscinet-getitem?mr=#1}{#2}
}
\providecommand{\href}[2]{#2}


\begin{thebibliography}{10}


\bibitem{Bing:Fixed}
\textit{{Bing, R. H.}}, {The elusive fixed point property}, {Amer. Math.
  Monthly} \textbf{{76}} ({1969}), {119--132}.

\bibitem{Bla:newk-theory}
\textit{{Blackadar, Bruce}}, {$K$-Theory for Operator Algebras},
  \textit{{Mathematical Sciences Research Institute Publications}}, volume~5,
  {Cambridge University Press}, 1998, {Second} edition.

\bibitem{BlaKumRor:apprdiv}
\textit{{Blackadar, Bruce and Kumjian, Alex and R\o rdam, Mikael}},
  {Approximately central matrix units and the structure of non-commutative
  tori}, {$K$-theory} \textbf{6} (1992), 267--284.

\bibitem{Brown:LefFixPt}
\textit{{Brown, Robert F.}}, {The Lefschetz Fixed Point Theorem}, {Scott,
  Foresman and Co.}, {1970}.

\bibitem{Con:class}
\textit{{Connes, Alain}}, {Classification of injective factors}, {Ann. of
  Math.} \textbf{104} (1976), 73--115.

\bibitem{Exel:ParAct}
\textit{{Exel, Ruy}}, {Circle actions on {$C^*$}-algebras, partial
  automorphisms, and a generalized {P}imsner-{V}oiculescu exact sequence}, {J.
  Funct. Anal.} \textbf{{122}} ({1994}), no.~{2}, {361--401}.

\bibitem{FatHer:Diffeo}
\textit{{Fathi, Albert and Herman, Michael}}, {Existence de
  diff{\'e}omorphismes minimaux}, {Ast{\'e}risque} \textbf{{49}} ({1977}),
  {37--59}.

\bibitem{GioKerr:Subshifts}
\textit{{Giol, Julien and Kerr, David}}, {Subshifts and perforation}, {J. Reine
  Angew. Math.} \textbf{{639}} ({2010}), {107--119}.

\bibitem{GioPutSkau:orbit}
\textit{{Giordano, Thierry and Putnam, Ian F. and Skau, Christian F.}},
  {Topological orbit equivalence and $C^*$-crossed products}, {J. Reine Angew.
  Math.} \textbf{{469}} ({1995}), {51--111}.

\bibitem{Hirsch:DiffTop}
\textit{{Hirsch, Morris}}, {Differential Topology}, \textit{{Graduate Texts in
  Mathematics}}, volume~{33}, {Springer-Verlag}, {1976}.

\bibitem{JiaSu:Z}
\textit{{Jiang, Xinhui and Su, Hongbing}}, {On a simple unital projectionless
  {$C^*$}-algebra}, {Amer. J. Math.} \textbf{{121}} ({1999}), no.~2,
  {359--413}.

\bibitem{KirWinter:dr}
\textit{{Kirchberg, Eberhard and Winter, Wilhelm}}, {Covering dimension and
  quasidiagonality}, {Internat. J. Math.} \textbf{15} ({2004}), no.~{1},
  {63--85}.

\bibitem{Lin:TAF1}
\textit{{Lin, Huaxin}}, {Tracially AF $C^{*}$-algebras}, Trans. Amer. Math.
  Soc. \textbf{353} (2001), 693--722.

\bibitem{Lin:asu-class}
\textit{{Lin, Huaxin}}, {Asymptotic unitary equivalence and classification of
  simple amenable {$C^*$}-algebras}, {Invent. Math.} \textbf{{183}} ({2011}),
  no.~{2}, {385--450}.

\bibitem{Lin:Spheres}
\textit{{Lin, Huaxin}}, {Minimal dynamical systems on connected odd dimensional
  spaces}, {2014}. {arXiv preprint math.OA/1404.7034}.

\bibitem{Lin:MinDyn}
\textit{{Lin, Huaxin}}, {Crossed products and minimal dynamical systems}, 2015.
  {arXiv preprint math.OA/1502.06658}.

\bibitem{LinNiu:KKlifting}
\textit{{Lin, Huaxin and Niu, Zhuang}}, Lifting {$KK$}-elements, asymptotic
  unitary equivalence and classification of simple {$C\sp \ast$}-algebras,
  {Adv. Math.} \textbf{{219}} ({2008}), no.~5, {1729--1769}.

\bibitem{LinNiu:RationallyTAI}
\textit{{Lin, Huaxin and Niu, Zhuang}}, {The range of a class of classifiable
  separable simple amenable C*-Algebras}, {J. Funct. Anal.} \textbf{{260}}
  ({2011}), no.~{1}, {1--29}.

\bibitem{LinPhi:MinHom}
\textit{{Lin, Huaxin and Phillips, N. Christopher}}, {Crossed products by
  minimal homeomorphisms}, {J. Reine Angew. Math.} \textbf{{641}} ({2010}),
  {95--122}.

\bibitem{LinQPhil:KthoeryMinHoms}
\textit{{Lin, Qing and Phillips, N. Christopher}}, {Ordered {$K$}-theory for
  {$C^\ast$}-algebras of minimal homeomorphisms}, in: {Operator algebras and
  operator theory ({S}hanghai, 1997)}, {Amer. Math. Soc.}, 1998,
  \textit{{Contemp. Math.}}, volume {228}, {289--314}.

\bibitem{MvN:Rings4}
\textit{{Murray, Francis J. and von Neumann, John}}, {On rings of operators.
  {IV}}, {Ann. of Math. (2)} \textbf{{44}} ({1943}), {716--808}.

\bibitem{Phillips:recsub}
\textit{{Phillips, N. Christopher}}, {Recursive subhomogeneous algebras},
  {Trans. Amer. Math. Soc.} \textbf{{359}} ({2007}), no.~{10}, {4595--4623
  (electronic)}.

\bibitem{Putnam:MinHomCantor}
\textit{{Putnam, Ian F.}}, {The $C^*$-algebras associated with minimal
  homeomorphisms of the Cantor set}, Pacific J. Math. \textbf{136} (1989),
  no.~2, 329--353.

\bibitem{Put:TopStabCantor}
\textit{{Putnam, Ian F.}}, {On the topological stable rank of certain
  transformation group {$C^*$}-algebras}, {Ergodic Theory Dynam. Systems}
  \textbf{{10}} ({1990}), no.~{1}, {197--207}.

\bibitem{Put:K-theoryGroupoids}
\textit{{Putnam, Ian F.}}, {On the {$K$}-theory of {$C^*$}-algebras of
  principal groupoids}, {Rocky Mountain J. Math.} \textbf{{28}} ({1998}),
  no.~{4}, {1483--1518}.

\bibitem{Ren:groupoid}
\textit{{Renault, Jean}}, {A groupoid approach to $C^*$-algebras},
  Springer-Verlag, 1980.

\bibitem{Ror:encyc}
\textit{{R{\o}rdam, Mikael}}, {Classification of nuclear, simple
  {$C^*$}-algebras}, in: {Classification of nuclear {$C^*$}-algebras. {E}ntropy
  in operator algebras}, {Springer}, {2002}, \textit{{Encyclopaedia Math.
  Sci.}}, volume {126}, {1--145}.

\bibitem{Ror:simple}
\textit{{R{\o}rdam, Mikael}}, {A simple $C^*$-algebra with a finite and an
  infinite projection}, {Acta Math.} \textbf{{191}} ({2003}), no.~{1},
  {109--142}.

\bibitem{SatWhiWin:NucDimZ}
\textit{{Sato, Yasuhiko and White, Stuart A. and Winter, Wilehlm}}, {Nuclear
  dimension and $\mathcal{Z}$-stability}, 2014. {arXiv preprint
  math.OA/1403.0747; to appear in Invent. Math.}

\bibitem{Str:XxSn}
\textit{{Strung, Karen R.}}, {On the classification of $C^*$-algebras of
  minimal product systems of the Cantor set and an odd dimensional sphere}, {J.
  Funct. Anal.} \textbf{268} (2015), no.~3, 671--689.

\bibitem{StrWin:Z-stab_min_dyn}
\textit{{Strung, Karen R. and Winter, Wilhelm}}, {Minimal dynamics and
  $\mathcal{Z}$-stable classification}, {Internat. J. Math.} \textbf{{22}}
  ({2011}), no.~{1}, {1--23}.

\bibitem{Toms:classproblem}
\textit{{Toms, Andrew S.}}, {On the classification problem for nuclear
  {$C^\ast$}-algebras}, {Ann. of Math. (2)} \textbf{{167}} ({2008}), no.~{3},
  {1029--1044}.

\bibitem{TomsWin:ssa}
\textit{{Toms, Andrew S. and Winter, Wilhelm}}, {Strongly self-absorbing
  {$C^*$}-algebras}, {Trans. Amer. Math. Soc.} \textbf{{359}} ({2007}),
  no.~{8}, {3999--4029}.

\bibitem{TomsWinter:PNAS}
\textit{{Toms, Andrew S. and Winter, Wilhelm}}, {Minimal dynamics and the
  classification of {$C^*$}-algebras}, {Proc. Natl. Acad. Sci. USA}
  \textbf{{106}} ({2009}), no.~{40}, {16942--16943}.

\bibitem{TomsWinter:minhom}
\textit{{Toms, Andrew S. and Winter, Wilhelm}}, {Minimal {D}ynamics and
  {K}-{T}heoretic {R}igidity: {E}lliott's {C}onjecture}, {Geom. Funct. Anal.}
  \textbf{{23}} (2013), no.~{1}, {467--481}.

\bibitem{Vil:perforation}
\textit{{Villadsen, Jesper}}, {Simple $C^*$-algebras with perforation}, {J.
  Funct. Anal.} \textbf{154} (1998), no.~1, {110--116}.

\bibitem{Wind:not_uniquely_ergo}
\textit{{Windsor, Alistair}}, {Minimal but not uniquely ergodic
  diffeomorphisms}, in: Smooth ergodic theory and its applications ({S}eattle,
  {WA}, 1999), Amer. Math. Soc., Providence, RI, 2001, \textit{Proc. Sympos.
  Pure Math.}, volume~69, 809--824.

\bibitem{Winter:lfdr}
\textit{{Winter, Wilhelm}}, {Simple $C^*$-algebras with locally finite
  decomposition rank}, {J. Funct. Anal.} \textbf{243} (2007), 394--425.

\bibitem{Winter:dr-Z-stable}
\textit{{Winter, Wilhelm}}, {Decomposition rank and $\mathcal{Z}$-stability},
  {Invent. Math.} \textbf{{179}} ({2010}), no.~{2}, {229--301}.

\bibitem{Win:Z-stabNucDim}
\textit{{Winter, Wilhelm}}, {Nuclear dimension and {$\mathcal{Z}$}-stability of
  pure {$C^{*}$}-algebras}, {Invent. Math.} \textbf{{187}} ({2012}), no.~{2},
  {259--342}.

\bibitem{Win:ClassCrossProd}
\textit{{Winter, Wilhelm}}, {Classifying crossed products}, 2013. {arXiv
  preprint math.OA/1308.5084}.

\end{thebibliography}
\end{document}